\documentclass[fleqn,reqno]{amsart}
\usepackage[utf8]{inputenc}
\usepackage{amsmath}
\usepackage{amsfonts}
\usepackage{amsthm}
\usepackage{amssymb}
\usepackage{eucal,mathrsfs}
\usepackage{enumerate}
\usepackage{url}
\usepackage{graphicx}
\usepackage{hyperref}
\usepackage[usenames,dvipsnames,svgnames,table]{xcolor}
\definecolor{badgerred}{rgb}{0.715,0.004,0.004} 
\definecolor{burntorange}{rgb}{0.801,0.332,0.0}
\usepackage[font=small,labelfont={bf,sf}, textfont={sf},margin=0em]{caption}
\hypersetup{colorlinks,
filecolor=black,
linkcolor=badgerred,
citecolor=badgerred,
urlcolor=RoyalBlue,
bookmarksopen=true}
\hypersetup{bookmarksopenlevel=1}

\theoremstyle{definition}

\numberwithin{equation}{subsection}

\newcommand{\N}{\mathbb{N}}
\newcommand{\R}{\mathbb{R}}
\newcommand{\Z}{\mathbb{Z}}

\newcommand{\vk}{{\boldsymbol{k}}}
\newcommand{\vT}{{\boldsymbol{T}}}

\newcommand{\lin}{{\mathscr L}}

\title{Which shapes can appear in a Curve Shortening Flow Singularity?}
\author[]{Sigurd B. Angenent}
\thanks{The work presented here was done in Spring and Fall of 2022 as part of the \textsc{Madison Experimental Mathematics} lab ( \url{https://mxm.math.wisc.edu/} ). We would like to thank Kaiyi Huang and Devanshi Merchant for their role in mentoring the two student groups that participated.}
\author[]{Evan Patrick Davis} 
\author[]{Ellie DeCleene} 
\author[]{Paige Ellingson} 
\author[]{Ziheng Feng} 
\author[]{Edgar Gevorgyan} 
\author[]{Aris Lemmenes} 
\author[]{Alex Moon} 
\author[]{Tyler Joseph Tommasi} 
\author[]{Yamin Zhou} 
\address{S.B.Angenent, E.DeCleene, P.Ellingson, Z.Feng, T.Tommasi, Y.Zhou:  University of Wisconsin--Madison}
\address{E.P.Davis: Department of Mathematics, University of California, Los Angeles}
\address{Edgar Gevorgyan: Department of Mathematics, Rice University }
\address{A.Lemmenes: Seattle }
\address{A.Moon: Department of Mathematical Sciences, University of Wisconsin--Milwaukee}

\date{\today}
\begin{document}

\begin{abstract}
We study possible tangles that can occur in singularities of solutions to plane Curve Shortening Flow.  We exhibit solutions in which more complicated tangles with more than one self-intersection disappear into a singular point.  It seems that there are many examples of this kind and that a complete classification presents a problem similar to the problem of classifying all knots in $\R^3$.  As a particular example, we introduce the so-called $n$-loop curves, which generalize Matt Grayson's Figure-Eight curve, and we conjecture a generalization of the Coiculescu-Schwarz asymptotic bow-tie result, namely, a vanishing $n$-loop, when rescaled anisotropically to fit a square bounding box, converges to a ``squeezed bow-tie,'' i.e.~the curve $\{(x, y) : |x|\leq 1, y=\pm x^{n-1}\}\cup\{(\pm 1, y) : |y|\leq 1\}$.  As evidence  in support of the conjecture, we provide  a formal asymptotic analysis on one hand, and a numerical simulation for the cases $n=3$ and $n=4$ on the other.
\end{abstract}

\maketitle

\section{Introduction}

\subsection{Definitions}
Curve Shortening Flow deals with evolving families of closed curves.  Such a family is given by a differentiable map $\gamma:[0, T)\times\R \to \R^2$ where $u\mapsto\gamma(t, u)$ parametrizes the curve at time $t\in[0, T)$.  We assume that the curves are closed, i.e.
\begin{equation} \label{eq:curve-closed} 
  \gamma(t, u+1)=\gamma(t, u) \text{ for all } t, u,
\end{equation}
and that the curves are \emph{immersed,} i.e.
\[
\gamma_u(t, u)\stackrel{\rm def}=\frac{\partial \gamma}{\partial u}\neq 0 \text{ for all } t, u.
\]
By definition, a family $\gamma(t,\cdot)$ of curves evolves by Curve Shortening Flow if it satisfies
\begin{equation}\label{eq:CSF}
\gamma_t(t, u)^\perp = \vk(t, u).
\end{equation}
Here $\gamma_t^\perp$ is the component of $\gamma_t$ that is perpendicular to $\gamma_u$, i.e.,
\[
\gamma_t^\perp = \gamma_t - \langle \gamma_t, \gamma_u \rangle \gamma_u/\|\gamma_u\|^2.
\]
The map $\gamma:[0, T)\times\R \to \R^2$ is called a \emph{normal parametrization} if $\gamma_t\perp \gamma_u$ for all $(t, u)$. For a normal parametrization it follows that the normal component of the velocity $\gamma_t^\perp$ coincides with the velocity $\gamma_t$ itself.  Therefore, a normal parametrization evolves by Curve Shortening Flow exactly whenever
\[
\gamma_t(t, u) = \vk(t, u)
\]
for all $(t, u)$.

\subsection{Existence and uniqueness}
It has been shown that for any initial curve $\gamma^0$ there is a unique solution to Curve Shortening Flow starting from $\gamma^0$. More precisely, for any given immersed continuously differentiable initial curve $\gamma^0:\R\to\R^2$, there exist a $T>0$ and a solution $\gamma:[0, T)\times\R \to \R^2$ of~\eqref{eq:CSF} with $\gamma(0, u) = \gamma^0(u)$ for all $u\in\R$.  The solution $\gamma$ consists of closed curves, i.e.,~$\gamma$ satisfies~\eqref{eq:curve-closed}.  The solution is also unique up to a re-parametrization.  This means that any other solution $\tilde\gamma:[0, T)\times\R \to \R^2$ is obtained by re-parameterizing $\gamma$ in the sense that there is a differentiable function $\phi:[0, T)\times\R \to \R$ with $\phi_u(t, u)>0$ and $\phi(t, u+1)=\phi(t,u)$ for all $t, u\in[0, T)\times\R$ such that
\begin{equation}
\label{eq:reparametrization}
\tilde\gamma(t, u) = \gamma(t, \phi(t, u)) \qquad \forall (t, u)\in[0, T)\times\R.
\end{equation}
The maximal time interval $[0,T)$ during which the solution $\gamma$ is defined depends on the initial curve, but it is always finite.  The solution becomes singular at $t=T$, in the sense that the maximal curvature $k_{\max}(t) = \max_{u\in\R}|k(t, u)|$ on the curve becomes infinite:
\begin{equation}
\label{eq:curvature-blowup}
\lim_{t\nearrow T}k_{\max}(t) = +\infty.
\end{equation}

\subsection{The Gage-Hamilton-Grayson theorem for embedded
curves}\label{sec:GHGthm} If the initial curve $\gamma(0,\cdot)$ is embedded,
i.e.~if $\gamma(0, \cdot)$ has no self-intersections, then its corresponding deformation by Curve Shortening Flow is fairly well understood.  Grayson~\cite{grayson:hes86} showed that if $\gamma(0, \cdot)$ is embedded, then $\gamma(t, \cdot)$ is embedded for all $t\in [0, T)$, and there is a time $t_c\in[0,T)$ such that $\gamma(t, \cdot)$ is convex for $t\in[t_c, T)$.  The Gage-Hamilton theorem~\cite{gagehamilton:hes86} guarantees that the curve $\gamma(t,\cdot)$ shrinks to a point $p\in\R^2$, and that the shape of $\gamma(t, \cdot)$ right before it shrinks to $p$ is that of a circle centered at $p$ with radius $\sqrt{2(T-t)}$.

\subsection{Oaks' Theorem for immersed curves}\label{sec:Oaks}
\begin{figure}[t]
\includegraphics[width=0.7\textwidth]{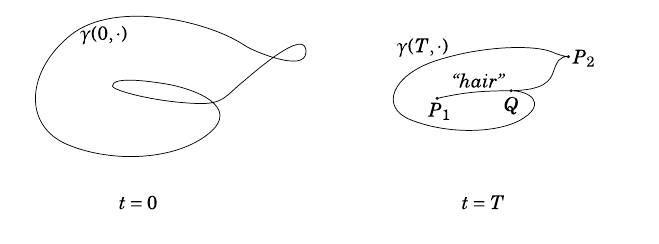}
  
\caption[]{According to Oaks' theorem \cite{Oaks1994} any solution $\gamma(t, \cdot)$ to Curve Shortening Flow that becomes singular at some finite time $T>0$ has a limit curve $\gamma(T, \cdot)$, which has finitely many singular points $P_1, \dots , P_N$. Right before the singularity happens, the curve $\gamma(t, \cdot)$ has a self-intersection arbitrarily close to each singular point.  In this figure above, the two loops in the initial curve contract at the same time, each resulting in one of the singular points $P_1, P_2$.  A curious possibility that Oaks' theorem leaves open is the existence of ``hairs'' such as the arc $P_1Q$ in the limit curve $\gamma(T, \cdot)$.}
\label{fig:Oaks}
\end{figure}

See Fig.~\ref{fig:Oaks}.  

Let $\gamma:[0, T)\times\R\to\R^2$ be a solution to \eqref{eq:CSF} whose initial curve is not necessarily embedded, and assume that $\gamma$ has been re-parameterized so that $\gamma_t\perp\vT$.  Then it was shown in \cite{angenent:pec90} that
\[
\gamma(T, u) = \lim_{t\nearrow T} \gamma(t, u)
\]
exists, and that the curves $\gamma(t, \cdot)$ remain smooth except near finitely many ``singular points'' $P_1, \dots, P_N\in\R^2$; more precisely, the curvature of $\gamma(t, \cdot) $ is uniformly bounded away from the singular points i.e.~for any $A>0$ there is an $\epsilon_A>0$ such that $|k(t, u)|\leq A$ if $\gamma(t, u)\not\in \cup_i B(P_i; \epsilon_A)$.  Oaks~\cite{Oaks1994} proved that for some sequence of times $t_k\nearrow T$, the curve $\gamma(t_k, \cdot)$ must have a self-intersection near each singular point $P_i$ ($1\leq i\leq N$).

\begin{figure}[t]
\centering \includegraphics[width=\textwidth]{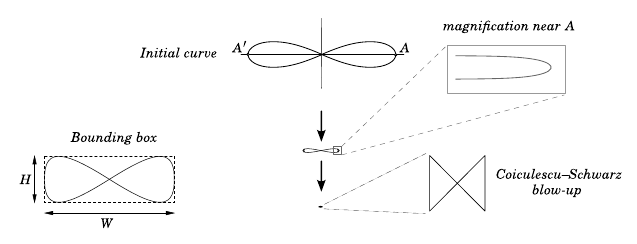}
\caption{Grayson's Symmetric Figure-Eight.}
\label{fig:SymmetricEight}
\end{figure}

\subsection{Two typical examples}
\begin{figure}[ht]
\centering \includegraphics[width=\textwidth]{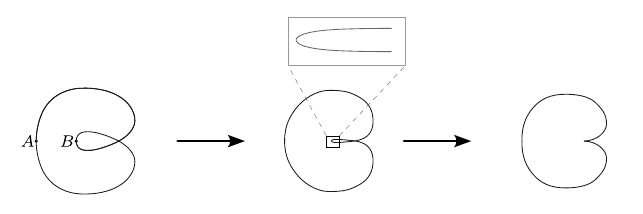}
\caption{The Symmetric Cardioid.}
\label{fig:Cardioid}
\end{figure}

There are two examples of curves with self-intersect\-ions whose evolution by Curve Shortening Flow has been studied in detail: the Symmetric Figure-Eight (Fig.~\ref{fig:SymmetricEight}), and the Cardioid (Fig.~\ref{fig:Cardioid}).  In both examples one self-intersection of the curve vanishes precisely into the singular point.

\subsubsection{The Symmetric Figure Eight}
See Fig.~\ref{fig:SymmetricEight}.  Matt~Grayson \cite{Grayson:figeight89} showed that a Figure-Eight curve with reflection symmetry in the origin, shrinks to a point under Curve Shortening Flow. He also showed that the curve becomes flat at the time of singularity $T$, in the sense that $\frac{H(t)}{W(t)} \to 0 $ as $t\nearrow T$, where $H(t)$ and $W(t)$ are height and width of the bounding box of the curve at time $t$. The bounding box is the smallest rectangle that encloses the curve.

If the curvature, at time $t=0$,  along the curve is strictly increasing between the point of self-intersection and the endpoint $A$ where the tangent is vertical, then this remains true for all $t\in[0, T)$.  Coiculescu and Schwartz \cite{coiculescu2022affine} recently showed for such curves that if one rescales the shrinking curve so that its bounding box becomes the unit square, then right before the singular time the resulting curve converges to a ``bow-tie.''

In a type-II blow-up one zooms in on the curve at the point $A$ of maximal curvature $k_{\max}(t)$ and magnifies the curve by a factor $k_{\max}(t)$.  In this case the magnified curve converges to the so-called ``Grim Reaper,'' a translating solution of Curve Shortening Flow, which is given by $x=\ln \cos y$.

\subsubsection{The Symmetric Cardioid}
See Fig.~\ref{fig:Cardioid}.  The Cardioid-like initial curve on the left has no inflection points, and only one self-intersection. It is symmetric with respect to reflection in the $x$-axis, and the curvature is strictly increasing along the curve between the points $B$ and $A$. It was shown in \cite{angenent:fsc91} that these properties are preserved by Curve Shortening Flow, and that the smaller inner loop contracts to a point before the outer loop can shrink too much.  The shape of the inner loop near its point of maximal curvature $B$ is again given by a Grim Reaper, as was the case for the symmetric Figure-Eight.  The precise rate at which the curvature at the point $B$ blows up is $k_{\max}(t) \sim \sqrt{\frac{\ln|\ln(T-t)|}{T-t}}$.  See~\cite{1995AngenentVelazquez}.

\begin{figure}[b]
\centering \includegraphics{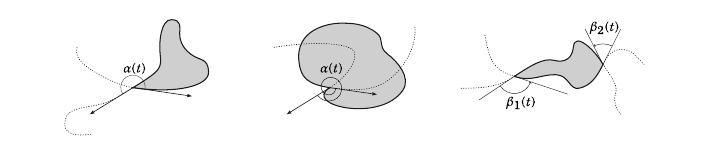}
\caption{A loop with a convex corner, one with a concave corner, and an eye}
\label{fig:areaofloop}
\end{figure}

\section{The rate at which loops and eyes lose area}
\subsection{Loops and eyes}
In both the Cardioid and Figure-Eight examples, exactly one self-intersection of the solution vanishes into the singularity.  Oaks' theorem~\ref{sec:Oaks} allows for the possibility that a singularity absorbs more than one self-intersection.  Our goal is to present and analyze a number of different examples of solutions where this does indeed occur.  Before describing our examples we need to discuss \emph{loops} and \emph{eyes} in plane immersed curves, and the rate at which their enclosed area decreases.

By definition, a \emph{loop} of a solution to Curve Shortening Flow.  $\gamma:(t_0, t_1)\times\R\to\R^2$ is a pair of functions $a, b:(t_0, t_1)\to\R$ such that
\begin{enumerate}[(i)]
\item $\gamma(t,a(t)) = \gamma(t, b(t))$ for all $t\in(t_0,t_1)$
\item $\gamma(t, \cdot)$ is injective on $[a(t), b(t)]$, i.e.~if $t\in(t_0, t_1)$ and $u, v\in[a(t), b(t)]$ then $\gamma(t, u)=\gamma(t, v)$ implies $u=v$.
\end{enumerate}
If $(a, b)$ is a loop for $\gamma$, then $u\in[a(t), b(t)]\mapsto \gamma(t, u)\in\R^2$ defines a closed curve without self-intersections.  The curve $\gamma(t, [a(t), b(t)])$ has a corner at $\gamma(t, a(t))$, where the two tangents $\gamma_u(t, a(t))$ and $\gamma_u(t, b(t))$ form an angle.  Depending on the tangents, the curve can have a \emph{convex corner} (Figure~\ref{fig:areaofloop} on the left) or a \emph{concave corner} (Figure~\ref{fig:areaofloop} in the middle).

\subsection{The area of a loop or eye}
If we let $A(t)$ be the enclosed area of the loop, then it is known that
\begin{equation}\label{eq:loop-area-loss}
\frac{dA}{dt} = -\int_{a(t)}^{b(t)} k(t, u) ds.
\end{equation}
In the case where the loop has a convex corner one has
\[
\int_{\text{loop}} kds = \alpha(t)
\]
where $\alpha(t)$ is the angle between the tangents $\gamma_u(t, a(t))$ and $\gamma_u(t, b(t))$ (see Figure~\ref{fig:areaofloop}). Since the angle $\alpha(t)$ always lies between $\pi$ and $2\pi$, we see that the area of a loop always satisfies
\begin{equation}\label{eq:loop-area-loss-bound}
-2\pi < \frac{dA}{dt} < -\pi.
\end{equation}
If the corner is concave (Figure~\ref{fig:areaofloop}, middle) then the same reasoning  shows again
\[
\frac{dA}{dt} = -\int_{\text{loop}}kds = -\alpha(t),
\]
where the angle $\alpha(t)$ now satisfies $2\pi<\alpha(t)<3\pi$.  For loops with a concave corner we therefore find
\[
-3\pi < \frac{dA}{dt} < -2\pi.
\]

One can similarly define \emph{an eye} of a solution to Curve Shortening Flow to be the union of two simple arcs in the evolving curve which share their two endpoints, but are otherwise disjoint.  An eye again defines a simple closed curve in the plane and thus encloses a certain area $A(t)$.  Just as for loops, the rate at which an eye loses area under Curve Shortening Flow is determined by the tangent angles at the corner points:
\begin{equation}
\label{eq:eye-area-loss}
\frac{dA}{dt} = -\int_{\text{eye}} k(t, u) ds = -2\pi+\beta_1(t)+\beta_2(t).
\end{equation}
The angles $\beta_i(t)$ are defined in Figure~\ref{fig:areaofloop}, and are constrained by $0<\beta_i(t)<\pi$.  This implies
\begin{equation}\label{eq:eye-area-loss-bound}
-2\pi < \frac{dA}{dt} <0.
\end{equation}

\section{Three intersections simultaneously vanishing into a singular point}
\label{sec:three-in-one}
\begin{figure}[h]
\includegraphics[width=\textwidth]{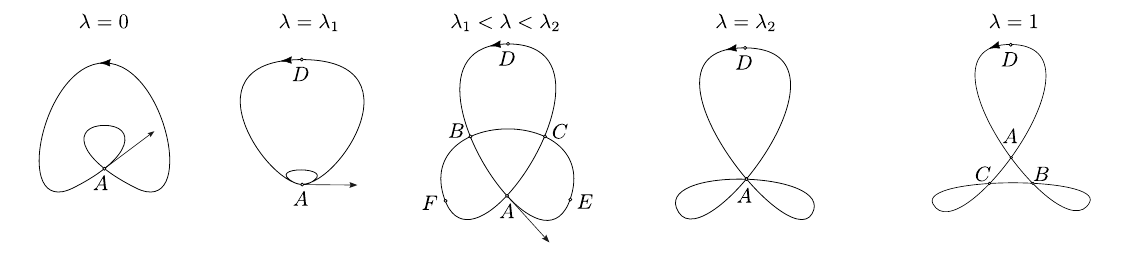}
\caption{Constructing a solution for which three self-intersections vanish into one singularity.}
\label{fig:three-intersections-vanish}
\end{figure}

We now show that a properly chosen initial curve will lead to a solution to Curve Shortening Flow that becomes singular without shrinking to a point, and for which three self intersections vanish into the singularity.  Our construction starts with a one parameter family of initial curves $\{\gamma_\lambda : 0\leq \lambda\leq 1\}$ (see Figure~\ref{fig:three-intersections-vanish}), and then uses an intermediate value argument to show that for at least one choice of the parameter $\lambda_*\in(0,1)$, the Curve Shortening Flow starting from $\gamma_{\lambda_*}$ does indeed have a singularity that absorbs three self-intersections.

The initial curves $\gamma_\lambda$ are immersed, locally convex curves with winding number~$2$ (their unit tangent rotates through $4\pi$ as one goes around the curves once). We also assume that they are symmetric with respect to reflection in the $y$-axis. As one increases $\lambda$ from $0$ to~$1$, the shape of the curve changes qualitatively twice, first at $\lambda=\lambda_1$, and then at $\lambda=\lambda_2$.  For $\lambda=\lambda_1$ the curve has a self-tangency (at $A$), and for $\lambda=\lambda_2$ the curve has a triple point (again labeled $A$).  For all other $\lambda\in[0,1]$, the curve has no self-tangencies or triple points.

Let $\gamma_\lambda(t)$ ($0\leq t <T_\lambda$) be the smooth solution to Curve Shortening Flow with $\gamma_\lambda$ as initial condition.  The life-time $T_\lambda$ of the solution depends on $\lambda$.

For $\lambda\in[0, \lambda_1)$, the initial curve has exactly one self-intersection, and since the number of self-intersections cannot increase under Curve Shortening Flow, the solution $\gamma_\lambda(t)$ has zero or one self-intersection for $t\geq 0$.  The winding number of $\gamma_\lambda(t)$ is $2$, so the curve is not embedded, and therefore has exactly one self-intersection for all $t\in[0, T_\lambda)$.

For $\lambda>\lambda_1$, the initial curve $\gamma_\lambda$ has three self intersections $A,B,C$, each of which is the convex corner of a loop in the curve (thus $A$ is the corner of the loop $ABDCA$; $B$ is the corner of the loop $BCEAB$; $C$ is the corner of the loop $CBFAC$.)  By symmetry, the $BCEAB$ and $CBFAC$ loops have equal area.  We assume that the area of the $ABCDA$ loop is strictly more than twice the area of the $BCEAB$ loop (and hence also more than twice the area of the $CBFAC$ loop). The rates at which all three loops lose area are bounded by~\eqref{eq:loop-area-loss-bound}, which implies that the smaller loops ($BCEAB$ and $CBFAC$) must vanish before the larger loop $ABDCA$ vanishes.

One way in which the smaller $CBFAC$ loop can vanish is if at some time $t_0>0$ the two intersections $A$ and $C$ come together and cancel.  When this happens, the arc $AEC$ becomes tangent to the arc $FACD$.  By symmetry, the intersection $B$ then also merges with $A$ and the arc $BFA$ has become tangent to $EABD$.  At this particular instant, i.e.,~at $t=t_0$, the curve looks like the initial curve $\gamma_\lambda$ with $\lambda=\lambda_1$. Immediately, after $t=t_0$ the solution will then have one self-intersection and will have the Cardioid shape of the initial curve with $\lambda=0$.  This leads us to consider the following set of parameter values
\[
I_1 \stackrel{\rm def}= \left\{ \lambda\in[0,1] \mid \exists t_1\in[0, T_\lambda) : \gamma_\lambda(t_1)\text{ has one transverse self-intersection}.  \right\}
\]
If an immersed curve has one transverse self-intersection, then a $C^1$ small perturbation of the curve also has exactly one transverse self-intersection.  By $C^1$-continuous dependence of solutions to Curve Shortening Flow on parameters, it follows that $I_1$ is an open subset of the interval $[0,1]$.

A second way in which the two smaller loops $CBFAC$ and $BAECB$ can vanish, occurs if the solution $\gamma_\lambda(t)$ develops a triple point at some time $t_0\in[0, T_\lambda)$.  This happens when the three intersection points $A, B, C$ come together (if $A$ and $B$ come together, then the reflection symmetry forces $C$ to meet $A$ and $B$ as well.)  Since the curves $\gamma_\lambda(t)$ are locally convex, the three branches of the curve passing through the triple point must keep moving as~$t$ increases, so that right after the triple point appears, the solution develops a triangle $ABC$ whose sides are concave (as in the initial curve with $\lambda=0$).  Under Curve Shortening Flow, such a concave triangle must keep growing. This prompts us to introduce
\[
I_3 \stackrel{\rm def}= \left\{ \lambda\in [0,1] \mid \exists t_1\in [0, T_\lambda) : \gamma_\lambda(t_1)\text{ has a concave triangle} \right\}.
\]
Continuous dependence of the solution on the parameter $\lambda$ implies that the set $I_3$ is also an open subset of $[0,1]$.

\begin{figure}[h]
\centering \includegraphics[width=0.5\textwidth]{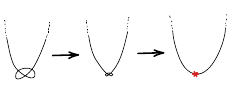}
\caption{Three intersections coming together in one singularity when $\lambda\in [0,1]\setminus (I_1\cup I_3)$. }
\label{fig:three-in-one}
\end{figure}

The two sets $I_1$ and $I_3$ are disjoint, and nonempty (because $0\in I_1$, $1\in I_3$). It follows that there exists at least one $\lambda_*\in[0,1] \setminus \left( I_1\cup I_3\right)$.  The solution $\gamma_{\lambda_*}(t)$ therefore has three self-intersections for all $t\in[0, T_\lambda)$, and never develops a triple point.  The latter implies that the intersection point $C$ never coincides with $A$, and that $C$ always lies between $A$ and $D$.  The solution $\gamma_{\lambda_*}$ thus always has a loop $CBFAC$, for all $t\in [0, T_{\lambda_*})$.  Since the area of this loop shrinks to zero as $t\nearrow T_{\lambda_*}$, then the loop must shrink to a point.  In particular, the intersections $A$ and $C$ converge to the same point as $t\nearrow T_{\lambda_*}$. By symmetry, the self-intersection $B$ also converges to the same point.

\section{Other singularity types}
One can imagine that there are more complicated figures that vanish into a singular point.  More specifically, consider a solution to Curve Shortening Flow $\bigl\{\gamma(t) \mid 0\leq t<T\bigr\}$ that forms a singularity at time $T$.  The number of self-intersections of $\gamma(t)$ does not increase with time, and thus it is eventually constant, i.e., ~there is an integer $n\in\N_0$ and $t_*\in (0, T)$ such that $\gamma(t)$ has exactly $n$ self-intersections for all $t\in(t_*, T)$.  One could then ask \emph{which curves with $n$ self-intersections can appear in a singular solution to Curve Shortening Flow?} 

To formulate more precise questions, we recall the notion of a \emph{flat knot} and a \emph{tangle}.  (Flat knots were defined in \cite{angenent:cst05}.)

\subsection{Flat knots and Tangles}
By definition, a flat knot in the plane is a closed immersed curve in $\R^2$ that has no self-tangencies (i.e.~all its self-intersections are transverse).  Two flat knots $\gamma_0, \gamma_1$ are defined to be equivalent if there is a smooth one parameter family $\{\gamma_\lambda \mid \lambda\in[0,1]\}$ of closed immersed curves such that each $\gamma_\lambda$ is a flat knot.  In other words, two flat knots are equivalent if one can be deformed into the other through immersed curves without ever forming a self-tangency.

Two equivalent flat knots have the same number of self-intersections, but the converse is not true: a Cardioid and a Figure-Eight both have one self-intersection, but they are not equivalent as flat knots.

We have several examples of plane Curve Shortening Flow in which the whole solution $\gamma(t)$ converges to one point. This happens for the circle, Grayson's symmetric Figure-Eight \cite{coiculescu2022affine,Grayson:figeight89}, as well as for the Abresch-Langer solitons \cite{abreschlanger:ncs86}. This suggests the following natural question:

\begin{quote}
Q1: \emph{Which other flat knots can shrink to a point under Curve Shortening Flow?}
\end{quote}

\begin{figure}[b]
\centering \includegraphics[width=0.7\textwidth]{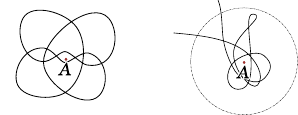}
\caption{Could this flat knot (left) or this tangle (right) vanish into a singularity at $A$?}
\label{fig:tangle}
\end{figure}

We also have examples of solutions that become singular without shrinking to a point.  To formulate the analogous questions, we define what a \emph{tangle} is.  Let $S\subset\R^2$ be a circle.  An immersed curve $\gamma$ \emph{has a tangle in $S$} if $\gamma$ intersects $S$ transversely in two distinct points, and if $\gamma$ has no self-tangencies inside $S$.  The tangles of two plane immersed curves $\gamma_0$, $\gamma_1$ in a circle $S$ are equivalent if there is a smooth one parameter family of plane immersed curves $\{\gamma_\lambda \mid \lambda\in[0,1]\}$ such that each $\gamma_\lambda$ has a tangle in $S$.

If a solution to Curve Shortening Flow does not shrink to a point, then Oaks' theorem provides us with a limit curve $\gamma_T=\lim_{t\nearrow T}\gamma_t$  with one or more singular points $P_1, \dots, P_N$, which are connected by smooth arcs.  For any $\epsilon>0$, we let $S_\epsilon(P_i)$ be the circle centered at $P_i$ with radius $\epsilon$.  If the solution has ``no hairs'' (see Section~\ref{sec:Oaks}, and Figure~\ref{fig:Oaks}), then there is an $\epsilon>0$ such that $\gamma_T$ intersects each circle $S_\epsilon(P_i)$ transversely in two points.  It follows that there is a $t_\dag\in(t_*, T)$ such that $\gamma(t)$ also intersects each $S_\epsilon(P_i)$ transversely in two points for $t\in(t_\dag,T)$.  Since $\gamma(t)$ has no self-tangencies it follows that $\gamma(t)$ has a tangle in each circle $S_\epsilon(P_i)$, and that for each $i$ the tangles of $\gamma(t)$ are equivalent for all $t\in (t_\dag, T)$.

For the Cardioid solution, the tangle that disappears in the singularity is that of a loop, and it has one self-intersection. For the solutions described in \S~\ref{sec:three-in-one}, the vanishing tangle is pretzel-shaped, and has three self-intersections. One can now ask the following question:
\begin{quote}
Q2: \emph{Which tangles can vanish into a singular point?}
\end{quote}
At present, we do not have answers to these questions and our goal here is to offer a number of examples of more complicated singularity types.

\section{$n$-loop curves}
\subsection{Definition of $n$-loops}
\begin{figure}[t]
\centering \includegraphics[width=0.8\textwidth]{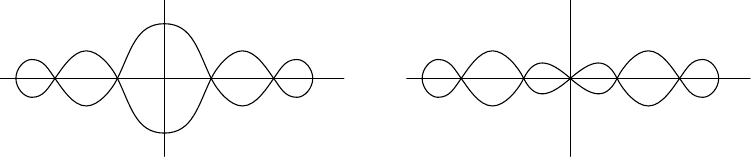}
\caption{$n$-loop curves with $n=5$ and $n=6$}
\label{fig:nLoops}
\end{figure}

Generalizing Figure-Eight curves, we consider $n$-loop curves (see Figure~\ref{fig:nLoops}), which we define to be closed immersed curves whose image is the union of two graphs
\[
y=\pm f(x)
\]
where $f:[-a, a]\to\R$ is a $C^1$ function that is either odd if $n$ is even, or else is even if $n$ is odd. The function $f$ has $n-1$ simple zeros in the interval $(-a,a)$. To ensure that the curve closes up smoothly (at least $C^1$), we also assume
\[
f(\pm a) = 0, \qquad \lim_{x\to a} f'(x) = -\infty.
\]
By symmetry, one then also has $\lim_{x\to -a}f'(x) = \pm\infty$, depending on the parity of $n$.

\subsection{The $n$-loop Conjecture}
\label{sec:nloop-conjecture}
For each $n\geq 3$, there exist $n$-loop curves whose evolution by Curve Shortening Flow consists of $n$-loop curves that shrink to the origin.

For $n=2$, an $n$-loop curve is a symmetric Figure-Eight, so the conjecture was already proved in Grayson's paper \cite{Grayson:figeight89}.  In fact, it follows from that paper that any solution starting from a symmetric Figure-Eight curve will remain a symmetric Figure-Eight while shrinking to the origin.  In contrast, it seems likely that most solutions to Curve Shortening Flow starting from $n$-loops with $n\geq 3$ will either lose self-intersections before they become singular, or else become singular without shrinking to a point.

\subsection{Proof of the $3$ and $4$-loop conjectures}
\label{sec:34loops-exist}
\begin{figure}[t]
\centering \includegraphics[width=\textwidth]{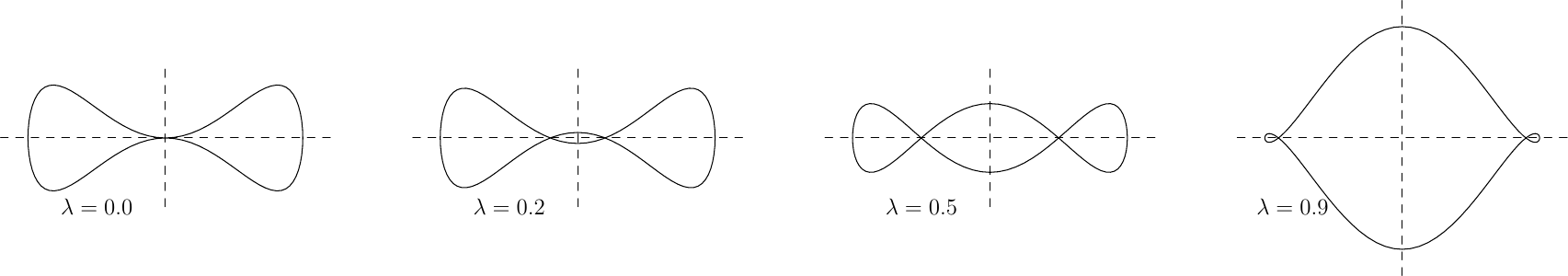}
\caption{Four of the initial curves $L_\lambda$ used in the construction of the vanishing three-loop.}
\label{fig:threeloop-proof}
\end{figure}

Consider the one parameter family of curves $L_\lambda$ ($\lambda\in[0,1)$) given by
\[
y^2 = (1-x^2)(x^2-\lambda^2)^2.
\]
One possible parametrization of these curves is
\[
L_\lambda(u) = \bigl(\cos u, (\cos^2 u - \lambda^2)\sin u\bigr).
\]
See Figure~\ref{fig:threeloop-proof}.

The curve $L_\lambda$ has two loops, one on the left, and the other on the right, and one ``eye'' in the middle.  It intersects the $x$-axis in the points $(\pm1, 0)$, which are the left and right extreme points on the curve, and also at the points $(\pm \lambda, 0)$, where $L_\lambda$ intersects itself transversely.

Let $\gamma_\lambda(t, u)$ be the solution to Curve Shortening Flow with $L_\lambda$ as initial condition.  Since Curve Shortening Flow is invariant under reflections, the solution $\gamma_\lambda$ is invariant under reflection in the $x$ and $y$-axes.  

The fact that the initial curve $L_\lambda$ is the union of two graphs is equivalent to the fact that the intersection of $L_\lambda$ with any vertical line $\{(x, y) : x=x_0\}$ contains at most two points.  The number of intersections of $\gamma_\lambda(t, \cdot)$ with any vertical line can not increase in time, so $\gamma_\lambda(t, \cdot)$ is always the union of two graphs
\begin{equation}
y=\pm u_\lambda(t, x), \qquad
-a_\lambda(t) < x < a_\lambda(t).
\end{equation}
The width $2a_\lambda(t)$ of the curve is a strictly decreasing function in time.

Since $L_\lambda$ has two self intersections, the solution $\gamma_\lambda(t, \cdot)$ has at most two self intersections.  Let us assume that $u_\lambda(t, x)>0$ for $x$ close to $a_\lambda(t)$.  In this case it follows that when $u_\lambda(t, 0)<0$ the function $x\mapsto u(t, x)$ has a zero $\bar x$ with $0 < \bar x < a_\lambda(t)$, so that the curve $\gamma_\lambda(t)$ is a $3$-loop.  If, on the other hand, $u(t, 0)>0$ then $x\mapsto u(t, x)$ has no zeros with $x<a_\lambda(t)$, and therefore the curve $\gamma_\lambda(t)$ is embedded.

\subsubsection{Claim} The set 
\[
I_1 \stackrel{\rm def}{=} \bigl\{\lambda \in [0, 1)\mid \exists t_0>0: \gamma_{\lambda}(t_0) \text{ is embedded}\bigr\}
\]
is nonempty and open.
\begin{proof}
Continuous dependence of the solution on initial data implies that $I_1$ is indeed open.

When $\lambda=0$, the initial curve $L_\lambda$ is embedded except for a self tangency at the origin, so the solution $\gamma_\lambda(t)$ immediately becomes embedded.  This implies $0\in I_1$, so $I_1$ is nonempty.
\end{proof}

\subsubsection{Claim} Let $A_{\lambda,\rm eye}(t)$ and $A_{\lambda,\rm ear}(t)$ denote the areas of the eye and the ear region of the solution $\gamma_\lambda(t)$.  Then the set
\[
I_2 \stackrel{\rm def}{=} \bigl\{
\lambda\in[0, 1] \mid \exists t_0\geq 0: A_{\lambda,\rm eye}(t_0) > 2 A_{\lambda,\rm ear}(t) >0
\bigr\}
\]
is nonempty and open in $[0,1)$.  

Note that the assumption $A_{\lambda,\rm ear}(t)>0$ implicitly assumes that $\gamma_\lambda(t)$ has self intersections.  Thus $I_1\cap I_2=\varnothing$.
\begin{proof}
Openness again follows from continuous dependence of the solution on initial data.  For $\lambda$ close to $\lambda=1$, the eye of the initial curve $L_\lambda$ has very small area, so there is a $\lambda_0\in(0, 1)$ such that $(\lambda_0,1)\subset I_2$ .
\end{proof}

We can now complete the proof of the $3$-loop conjecture.
Since $I_1$ and $I_2$ are disjoint nonempty open subsets of $[0, 1)$, there is at least one $\lambda_1\in(0, 1)$ that belongs to neither set. Now, consider the solution $\gamma_{\lambda_1}(t)$, which is defined for $t\in[0, T_{\lambda_1})$.   Since $\lambda_1\not\in I_{1}$ the curve $\gamma_{\lambda_1}(t)$ is a $3$-loop for all $t\in [0, T_1)$, so that it has an eye and ears for all $t$.  Furthermore, $\lambda_1\not\in I_2$ implies $A_{\lambda_1,\rm eye}(t)\leq A_{\lambda_1, \rm ear}(t)$ for all~$t\in [0, T_{\lambda_1})$.

The solution $\gamma_{\lambda_1}(t,\cdot)$ is defined for $t\in[0, T_{\lambda_1})$, and becomes singular as $t\nearrow T_{\lambda_1}$.   
It was shown in \cite{angenent:pec90,angenent:pec91} that the limit curve $\gamma_{\lambda_1}(T_{\lambda_1}, u) = \lim_{t\nearrow T_{\lambda_1}}(t, u)$ exists, and that it is smooth except at finitely many points.  In particular  $a_{\lambda_1}(T_{\lambda_1}) = \lim_{t\nearrow T_{\lambda_1}} a_{\lambda_1}(t)$ exists, and if $a_{\lambda_1}(T_{\lambda_1})>0$ then the limit curve is the union of two graphs $y=\pm u_{\lambda_1}(T_{\lambda_1}, x)$ ($|x|< a_{\lambda_1}(T_{\lambda_1})$).

By Oaks' theorem \ref{sec:Oaks}, the ears of the curve $\gamma_{\lambda_1}(t)$ must contract, and their area $A_{\lambda_1,\rm ear}(t)$ must go to zero.  Since $\lambda_1\not\in I_2$ the area of the eye also must vanish as $t\nearrow T_{\lambda_1}$.  If $a_{\lambda_1}(T_{\lambda_1}) > 0$ then this implies that $u_{\lambda_1}(T_{\lambda_1}, x) = 0$ and the limiting curve is the line segment between the two points $(\pm a_{\lambda_1}(T_{\lambda_1}), 0)$.  If this happens, then the limit curve $\gamma_{\lambda_1}(T_{\lambda_1})$ has a ``hair'', as in Figure~\ref{fig:Oaks}; the hair being the straight line segment connecting~$(\pm a_{\lambda_1}(T_{\lambda_1}), 0)$.  Grayson \cite{grayson:dinkytown89} has shown that this cannot happen, so we conclude that the solution $\gamma_{\lambda_1}(t, u)$ shrinks to the origin as $t\nearrow T_{\lambda_1}$.  This completes the proof of the $3$-loop conjecture.

\begin{figure}[h]
\centering \includegraphics[width=\textwidth]{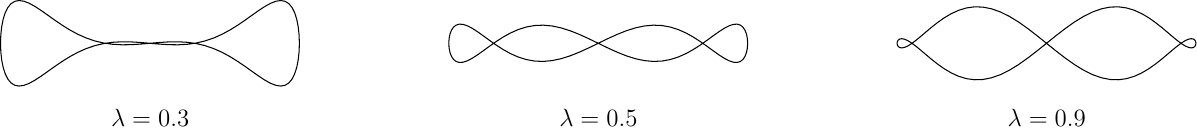}
\caption{Three of the initial curves $M_\lambda$ used in the construction of the vanishing four-loop.}
\label{fig:fourloop-proof}
\end{figure}

The $4$-loop conjecture can be proved in the same way by considering a one parameter family of initial curves depicted in Figure~\ref{fig:fourloop-proof}, and given by
\[
y^2 = x^2(1-x^2)(x^2-\lambda^2).
\]
These curves have two eyes and two ears.  Arguing as in the $3$-loop case, we find an exceptional value of $\lambda$ for which $\gamma_\lambda(t)$ has three self intersections at all times, and for which the area of either eye is no more than twice the area of either loop.  For such a $\lambda$, the solution remains a $4$-loop at all times, while it must also shrink to a point .

\section{$n$-loop curves--formal analysis}
\subsection{Nearly flat $n$-loops}
We derive a formal approximation to the solutions we have computed.  We begin by recalling that an $n$-loop curve that is symmetric with respect to reflection in both $x$ and $y$ axes is the union of two graphs
\begin{equation}
y=\pm u(t, x), \qquad \bigl(-a(t)\leq x\leq a(t),\quad 0\leq t <T\bigr).
\end{equation}
The left and right endpoints of the $n$-loop curve at time $t$ are $(\pm a(t), 0)$.  The curve evolves by Curve Shortening Flow if the function $u$ satisfies
\begin{equation}
u_t = \frac{u_{xx}}{1+u_x^2}, \qquad\bigl( -a(t)<x<a(t),\quad 0\leq t < T\bigr)
\end{equation}
At the end points $x=\pm a(t)$ the curve must have a vertical tangent, which leads to the boundary condition
\begin{equation}
u(t, a(t)) = \pm \infty,  \qquad u(t, -a(t)) = \pm\infty.
\end{equation}
Depending on whether $n$ is even or odd, the function $x\mapsto u(t, x)$ is odd or even.

\begin{figure}[h]
\centering \includegraphics{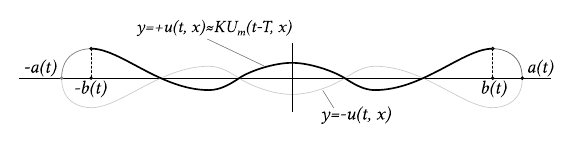}
\caption{$a(t)$ and $b(t)$ in a nearly five-loop.}
\label{fig:a-and-b}
\end{figure}
If we include the endpoints, then for an $n$-loop curve the function $x\mapsto u(t, x)$ has $n+1$ zeros.  On each loop the function $x\mapsto u(t, x)$ attains a maximum or minimum (depending on the sign of $u(t, x)$), so that the number of zeros of $x\mapsto u_x(t, x)$ is at least $n$.  We only consider initial curves where $x\mapsto u_x(0, x)$ has precisely $n$ zeros, so, since this number of zeros cannot increase in time, there is exactly one zero of $x\mapsto u_x(t, x)$ on each loop.

Let the local maximum/minimum of $x\mapsto u(t, x)$ on the right ear be at $x=b(t)$, so $u_x(t, b(t)) = 0$.  The local maximum/minimum on the left ear then occurs at $x=-b(t)$.

To arrive at our approximation, we now assume that $u_x$ is small on the part of the curve on which $|x|\leq b(t)$.  We can then approximate $1+u_x^2\approx 1$, and we find that for $|x|\leq b(t)$ the function $u$ should satisfy
\begin{equation}\label{eq:1}
u_t = \frac{u_{xx}}{1+u_x^2} \approx u_{xx}.
\end{equation}

We now recall some special solutions of the linear heat equation, which, as we will suggest, should describe the appearance of the $n$-loop solutions close to their singular time.  One finds these solutions by looking for polynomial solutions of the form
\[
u(t, x) = x^m + c_1(t)x^{m-1}+c_2(t)x^{m-2}+\cdots+c_{m-1}(t)x+c_m(t)
\]
where $c_1(t), \dots, c_m(t)$ are functions of time which are to be determined.  The first few functions of this type that appear are
\begin{align*}
  U_0(t, x) &= 1 &
                   U_1(t, x) &= x \\
  U_2(t, x) &= x^2+2t &
                        U_3(t, x)&=x^3+6xt\\
  U_4(t, x)&=x^4+12x^2t+12t^2 &
                                U_5(t, x) &=x^5+20x^3t+60t^2
\end{align*}
More generally one has:
\subsection{Lemma}
For each integer $m\geq 0$, the function
\begin{align}\label{eq:4}
  U_m(t, x)
  &=\sum_{k=0}^{\lfloor m/2\rfloor} \frac{m(m-1)\cdots(m-2k+1)}{k!} x^{m-2k}t^k\\
  &=x^m + \frac{m(m-1)}{1}x^{m-2}t
    + \frac{m(m-1)\cdots (m-3)}{1\cdot 2} x^{m-4}t^2+\cdots
    \notag
\end{align}
is a solution of the heat equation.
  
For $t>0$, the function $x\mapsto U_m(t, x)$ has no zeros, except $x=0$ in the case that $m$ is odd.
  
For $t<0$, the function $x\mapsto U_m(t, x)$ has $m$ zeros.
  
For any $K, T\in \R$ the function $\tilde U(t, x) = K U_m(t-T, x)$ is also a solution of the heat equation.

\begin{proof}
By direct substitution, one verifies $ \bigl(U_m\bigr)_t = \bigl(U_m\bigr)_{xx}$, i.e., ~that $U_m$ does indeed satisfy the heat equation. Another direct substitution will also show that $\tilde U$ satisfies the heat equation.
  
Assume $t>0$.  Then all coefficients in $U(t, c)$ are positive, so that $U_m(t, x)>0$ for all $x>0$.  If $m$ is even then $U_m(t, 0)>0$ also holds, and, since $x\mapsto U_m(t, x)$ is even, we also have $U_m(t, x)>0$ for all $x<0$.  If $m$ is odd, then $U_m(t, x) = -U_m(t, -x)<0$ for all $x>0$, while $U_m(t, 0)=0$.  We see that for odd $m$ the only zero of $x\mapsto U_m(t, x)$ is $x=0$.
  
For $t<0$, one can write the solutions $U_m$ in terms of Hermite polynomials $H_m(x)$ (see \cite{enwiki:1185297262}):
\begin{equation}\label{eq:6}
U_m(t, x) =  (-t)^{m/2}H_m\left(\frac{x}{2\sqrt{-t}}\right)
\end{equation}
The Hermite polynomials are orthogonal polynomials, which implies that $H_m$ has precisely $m$ zeros.
\end{proof}

\subsection{Conjecture concerning the shape of $n$-loops}\label{sec:conj-conc-shape}
Consider one of the $n$-loop solutions that shrink to the origin, and whose existence is conjectured in \S\ref{sec:nloop-conjecture}.

Then, when $t$ is close to $T$, the solution $u$ in the region $|x|\leq b(t)$, satisfies
\begin{equation}
  \label{co:1}
  u(t, x) \approx K U_{n-1}(t-T, x),
\end{equation}
where $K$ is a constant that depends on which $n$-loop solution we are considering.  Furthermore, $a(t)$ and $b(t)$ satisfy
\begin{align}
  &\lim_{t\nearrow T}\dfrac{b(t)}{a(t)} =1,   \label{co:2}\\
  &a(t) = (1+o(1)) \left(\frac{n\pi}{2K}(T-t)\right)^{1/n}.  \label{co:3}
\end{align}
Finally, the Coiculescu--Schwartz rescaling of the solution is
\begin{equation}
  \lim_{t\nearrow T} \frac{u(t, \xi a(t))}{|u(t, b(t))|} = \xi^{n-1}. \label{co:4}
\end{equation}
Since the bounding box of any solution that fits the description in Figure~\ref{fig:a-and-b} is given by $|x|\leq a(t)$, $|y|\leq u(t, b(t))$, this last statement~\eqref{co:4} says that in the limit $t\nearrow T$, the figure that appears upon rescaling the curve anisotropically such that its bounding box becomes the square $[-1,1]\times[-1,1]$, are the graphs of $y=\pm x^{n-1}$. This generalizes the Coiculescu--Schwartz theorem in \cite{coiculescu2022affine}, which treats the case $n=2$ (a Figure-Eight is a $2$-loop).

\subsection{Motivation for the conjecture}
One way to arrive at the approximation~\eqref{co:1} is to rewrite the partial differential equation~\eqref{eq:1} in self-similar variables
\[
X(t, x) = \frac{x}{\sqrt{T-t}}, \qquad U=\frac{u}{\sqrt{T-t}}, \qquad \tau(t)=-\ln (T-t).
\]
Note that as one approaches the singular time one has $t\nearrow T$ in the original variables, but $\tau\to\infty$ in the similarity variables. To compute the equation for $U(\tau, X)$, one substitutes
\[
\frac{u(t, x)}{\sqrt{T-t}} = U\bigl(\tau(t), X(t, x)\bigr) = U\Bigl(-\ln(T-t), \frac{x}{\sqrt{T-t}}\Bigr)
\]
in \eqref{eq:1} and applies the chain rule several times.  The result is
\begin{equation}
\label{eq:2}
U_\tau = \frac{U_{XX}}{1+U_X^2} - \frac{X}{2} U_X + \frac 12 U.
\end{equation}
Since
\begin{align*}
  u_x(t, x) &= \frac{\partial}{\partial x}\left(\sqrt{T-t} U(-\ln(T-t), x/\sqrt{T-t})\right)\\
            &=U_X\left(\sqrt{T-t} U(-\ln(T-t), x/\sqrt{T-t})\right)\\
            &=U_X(\tau, X)
\end{align*}
the assumption that $u_x$ be small for $|x|\leq b(t)$ is equivalent to the assumption that $|U_X|\ll 1$.  We can thus approximate the equation~\eqref{eq:2} by
\begin{equation}
\label{eq:3}
U_\tau \approx U_{XX} - \frac{X}{2} U_X + \frac 12 U.
\end{equation}
The right-hand side in this equation defines a linear operator
\[
\lin(U) \stackrel{\rm def}= U_{XX}-\frac X2 U_X + \frac 12 U,
\]
whose eigenfunctions are exactly the Hermite polynomials $H_m(X/2)$, with $\lambda_m = -(m-1)/2$ as corresponding eigenvalues. Assuming that $U(\tau, X)$ is described by solutions of the linear equation $U_\tau = \lin(U)$, we find
\[
U(\tau, X) = \sum_{m=0}^\infty c_m e^{-\lambda_m\tau}H_m(X/2),
\]
where $c_m\in\R$ are arbitrary coefficients, at least one of which is not zero.  Let $m_0\geq 0$ be the smallest value of $m$ for which $c_m\neq 0$.  We then have
\[
U(\tau , X) = \sum_{m=m_0}^\infty c_{m} e^{\lambda_{m}\tau} H_{m}(X/2).
\]
In the limit $\tau\to\infty$, the first term in this series is the largest.  Indeed,
\begin{multline*}
\sum_{m=m_0}^\infty c_{m} e^{\lambda_{m}\tau} H_{m}(X/2)
= e^{\lambda_{m_0}\tau} \sum_{m\geq m_0}c_{m} e^{(\lambda_{m}-\lambda_{m_0})\tau} H_{m}(X/2) \\
= e^{\lambda_{m_0}\tau} \sum_{m\geq m_0}c_{m} e^{-({m}-{m_0})\tau/2} H_{m}(X/2) \quad(\text{use }\lambda_m = -(m-1)/2).
\end{multline*}
For $m>m_0$, one has $\lim_{\tau\to\infty} e^{-(m-m_0)\tau} = 0$, so the only term in the series that does not vanish when $\tau\to\infty$ is the first term.  We therefore approximate our solution $U(\tau, X)$ by
\[
U(\tau, X) \approx c_{m_0}e^{-(m_0-1)\tau/2} H_{m_0}(X/2).
\]
After converting this approximation to the original variables $t, x,$ and $u$, one finds \eqref{co:1}, where $m_0=n-1$.

The condition~\eqref{co:2} is inspired by the analogous statements for the Cardioid in~\cite{1995AngenentVelazquez} and Figure-Eight in~\cite{coiculescu2022affine}: based on these we expect the outer loops, i.e., ~ the ones on the left and on the right that end at $x=\pm a(t)$ to have most of the area, and to be capped by translating solitons, i.e., ~by Grim Reapers.

To plot the Coiculescu--Schwartz rescaling of the $n$-loop, we note that the largest values of $x$ and $u$ on the curve are $a(t)$ and $|u(t, b(t))|$.  Then we plot
\[
u_{\rm CS} = \frac{u(t, x)}{|u(t, b(t))|} \text{ against } \xi = \frac{x}{a(t)}.
\]
Both $u_{\rm CS}$ and $\xi$ satisfy $|u_{\rm CS}|\leq 1$ and $|\xi|\leq 1$. Assuming the approximation \eqref{co:1} holds, we find
\begin{align}
  u_{\rm CS} &\approx \frac{KU_{n-1}(t-T, \xi a(t))}{KU_{n-1}(t-T, b(t))} \notag \\
             &=\frac{(\xi a(t))^{n-1} + c_{1}(\xi a(t))^{n-3}(t-T)  + c_{2}(\xi a(t))^{n-5}(t-T)^2 + \cdots}
               {b(t)^{n-1} + c_{1}b(t)^{n-3}(t-T) + c_{2}b(t)^{n-5}(t-T)^2 + \cdots}\label{eq:5}
\end{align}
where $c_{k} = \frac{(n-1)(n-2)\cdots(n-2k)}{k!}$ are the coefficients that define the polynomial $U_{n-1}$ (see~\eqref{eq:4}). Divide the numerator and denominator in \eqref{eq:5} by $b(t)^{n-1}$ and use
\[
\frac{a(t)}{b(t)} \approx 1, \text{ and } b(t) \gg \sqrt{T-t} \text{ whenever }t\approx T.
\]
Then we obtain
\[
u_{\rm CS}\approx \frac{\xi^{n-1}\bigl(\frac{a(t)}{b(t)}\bigr)^{n-1} + c_{ 1}\xi^{n-3}\bigl(\frac{a(t)}{b(t)}\bigr)^{n-3}\frac{t-T}{b(t)^2}+\cdots} {1+ c_{1}\frac{t-T}{b(t)^2} + \cdots} \approx \xi^{n-1}
\]
This is why we conjecture that \eqref{co:4} holds.

We conclude by explaining why we believe the size of the $n$-loop should be given by~\eqref{co:3}. The rightmost loop in the $n$-loop has its corner at the largest zero of $x\mapsto u(t, x)$.  Let $x_*(t)$ be the $x$-coordinate of this rightmost zero.  In view of the approximation $u(t, x)\approx KU_m(t-T, x)$ and \eqref{eq:6}, which expresses $U_m$ in terms of the $m^{\rm th}$ Hermite polynomials, we have
\begin{equation}\label{eq:7}
\frac{x_*(t)}{2\sqrt{T-t}}\approx z_{n-1},
\end{equation}
where
\[
z_m\stackrel{\rm def}= \max\{z\in\R | H_m(z)=0\}
\]
is the largest zero of $H_m$.

The rate at which the area of the rightmost loop is decreasing is  $\alpha(t)$.

\begin{figure}
\centering
\includegraphics{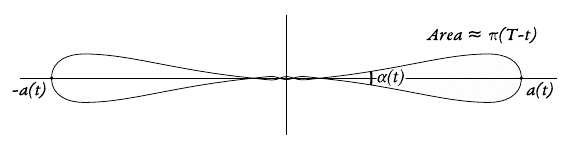}
\caption{\textbf{A $5$-loop just before the singularity. }  The shaded region is the right eye, whose area decreases at the rate $\pi-\alpha(t) = \pi-o(1)$, which implies that its area is $(\pi+o(1))(T-t)$.}
\label{fig:CS-n-loop}
\end{figure}

(See Figures~\ref{fig:areaofloop} and \ref{fig:CS-n-loop}).  In our case,
\[
\alpha(t) = \pi + 2\arctan u_x(t, x_*(t)).
\]
Using again the approximation $u\approx KU_m(t-T, x) $ , we find
\begin{align*}
    u_x(t-T, x_*(t))
   &\approx K\frac{\partial U_m}{\partial x}(t-T, x_*(t))\\
  &=\sum_{k=0}^{\lfloor m/2\rfloor} \frac{m(m-1)\cdots(m-2k)}{k!} x^{m-2k-1}(t-T)^k.
\end{align*}
It follows from \eqref{eq:7} that $x_*(t)\approx 0$ when $t\approx T$, which, in turn, implies $u_x(t-T, x_*(t))\approx 0$ and thus, $\alpha(t) \approx \pi$ for $t$ close to $T$.  The area of the rightmost loop is therefore
\begin{equation}\label{eq:8}
A_{\rm loop}(t) \approx \pi(T-t)
\end{equation}
One can also approximate the area enclosed by the right outer loop by computing the area it occupies in the Coiculescu--Schwartz rescaling, and then multiplying it with the scaling factors.
The region occupied by the rightmost loop in the Coiculescu--Schwartz scaling is approximately
\[
\bigl\{(\xi, u_{\rm CS}) : 0\leq \xi\leq 1, |u_{\rm CS}| \leq \xi^{n-1}\bigr\}
\]
whose area is
\[
\mathop{\mathrm{Area}} = \int_0^1 2\xi^{n-1}d\xi = \frac{2}{n}.
\]
Rescaling gives us the area of the rightmost loop in the $n$-loop, namely,
\[
A_{\rm loop}(t) \approx \frac 2n \times a(t) \times u(t, b(t))
\]
The approximation of $u(t, x)$ by $KU_{n-1}(t, x)$ tells us
\begin{align*}
  u(t, b(t))&\approx K\left\{ b(t)^{n-1} + c_{1}b(t)^{n-3}(t-T) + c_{2}b(t)^{n-5}(t-T)^2 + \cdots\right\}\\
            &=Kb(t)^{n-1}\left\{1 + c_1 \frac{t-T}{b(t)^2} + c_2\frac{(t-T)^2}{b(t)^4} + \cdots\right\}\\
            &\approx K b(t)^{n-1}
\end{align*}
because $b(t)\gg \sqrt{T-t}$ when $t\approx T$.  Using $b(t)\approx a(t)$, we then obtain
\[
A_{\rm loop}(t) \approx \frac{2}{n} K a(t)\times a(t)^{n-1} = \frac{2K}{n} a(t)^n.
\]
Combined with the other estimate~\eqref{eq:8} for the area, we conclude with
\[
\frac{2K}{n}a(t)^n \approx \pi (T-t)
\]
 and implies~\eqref{co:3}.
\section{$n$-loop curves--numerical study}
We ran numerical simulations of Curve Shortening Flow with the goal to find $3$-loop and $4$-loop solutions that shrink to a point (whose existence we proved in \S\ref{sec:34loops-exist}), and to verify the claims in \S\ref{sec:conj-conc-shape}.

\subsection{The numerical scheme}\label{sec:numerical-scheme}
To approximate parametrized solutions $\gamma:[0, T)\times\R \to\R^2$ to Curve Shortening Flow numerically, we start with the fact that any solution to
\begin{equation}\label{eq:deTurck}
\frac{\partial\gamma }{\partial t} = \frac{\gamma_{uu}}{\|\gamma_u\|^2}
\end{equation}
also satisfies Curve Shortening Flow \eqref{eq:CSF}.  Indeed, the curvature of the parametrized curve is
\[
\vk = \frac{1}{\|\gamma_u\|}\left(\frac{\gamma_u}{\|\gamma_u\|}\right)_u =\frac{\gamma_{uu}}{\|\gamma_u\|^2} - \left(\frac{1}{\|\gamma_u\|}\right)_u \frac{\gamma_u}{\|\gamma_u\|}.
\]
The second term is tangential to the curve, so~\eqref{eq:deTurck} implies $\left(\gamma_t\right)^\perp = \vk$, i.e., ~\eqref{eq:CSF}.  It was DeTurck's observation\footnote{In \cite{MR0697987}, DeTurck showed how the Ricci Flow can be reduced to a strictly parabolic equation.  The fact that \eqref{eq:deTurck} is equivalent to Curve Shortening Flow is a simpler version of DeTurck's reasoning in \cite{MR0697987}, and is generally known as ``DeTurck's trick.''}  that \eqref{eq:deTurck}~is a strictly parabolic equation, which can be used to prove short time existence for Curve Shortening Flow. Here we discretize \eqref{eq:deTurck} to approximate solutions to Curve Shortening Flow. In our numerical computation, we keep track of the values
\[
\gamma_i(t) = \gamma(t, \frac{i}{N}), \quad (i\in\Z)
\]
at a discrete set of times
\[
t=t_0, t_1, t_2, \dots.
\]
At each time $t_i$, we have a discrete approximation to the curve $\gamma(t_i, \cdot)$ and we compute an approximation at the next time $t_{i+1}$ by using a Backward Euler scheme.  In \eqref{eq:deTurck}, we approximate the second derivative by a second difference:
\[
\gamma_{uu}(t_{j+1}, i/N) \approx N^2 \bigl\{\gamma_{i+1}(t_{j+1}) - 2\gamma_i(t_{j+1}) + \gamma_{i-1}(t_{j+1})\bigr\}
\]
To approximate $\|\gamma_u(t, u)\|^2$, we consider the left and right differences
\[
D_r\gamma_i(t_j) = N\bigl\{\gamma_{i+1}(t_j) - \gamma_{i}(t_j)\bigr\}, \qquad D_l\gamma_i(t_j) = N\bigl\{\gamma_{i}(t_j) - \gamma_{i-1}(t_j)\bigr\}
\]
and set
\[
\|\gamma_u(t, u)\|^2 \approx \frac{1}{2}\left\{\|D_r\gamma_i(t_j)\|^2 + \|D_l\gamma_i(t_j)\|^2\right\}.
\]
Finally, we approximate the time derivative simply by
\[
\gamma_t(t_j, i/N) \approx \frac{\gamma_i(t_{j+1}) - \gamma_i(t_j)}{\Delta t_j}, \qquad \text{with }\Delta t_j = t_{j+1} - t_j.
\]
After replacing the partial derivatives in \eqref{eq:deTurck} with their discrete approximations, we arrive at a finite difference equation for the approximations $\gamma_i(t_j)$:
\begin{equation}
\label{eq:finite-diff}
\frac{\gamma_i(t_{j+1}) - \gamma_i(t_j)}{\Delta t_j}
=
2\,\frac{\gamma_{i+1}(t_{j+1}) - 2\gamma_i(t_{j+1}) + \gamma_{i-1}(t_{j+1})}
{\|\gamma_{i+1}(t_j)- \gamma_i(t_j)\|^2 + \|\gamma_{i}(t_j)- \gamma_{i-1}(t_j)\|^2}
\end{equation}
Note that we have somewhat arbitrarily decided to use the values of $\gamma$ at time $t_j$ to approximate the quantities $\|\gamma_u\|^2$, which appear in the denominator in \eqref{eq:deTurck}, while we approximated the second derivative $\gamma_{uu}$ using values of $\gamma$ at the next time $t_{j+1}$.  The advantage of doing this is that the finite difference equations \eqref{eq:finite-diff} can be written as a linear system for the unknown values of $\gamma$ at time $t_{j+1}$.  Namely, if we abbreviate
\[
K_i^j = \frac{2\Delta t_j}{\|\gamma_{i+1}(t_j)- \gamma_i(t_j)\|^2 + \|\gamma_{i}(t_j)- \gamma_{i-1}(t_j)\|^2} \quad\text{ and }\quad \gamma_i^j = \gamma_i(t_{j}),
\]
then \eqref{eq:deTurck} is equivalent to
\begin{equation}
\label{eq:linear-system}
-K_i^j\gamma_{i-1}^{j+1}+(1+2K_i^j)\gamma_i^{j+1}  - K_i^j\gamma_{i+1}^{j+1}
= \gamma_i^j.
\end{equation}
Keeping in mind that $\gamma_{i+N}(t) = \gamma_i(t)$, this is a finite system of linear equations whose matrix form is
\[\small
\begin{bmatrix}
1+2K_1^j & - K_1^j &0 &  &\dots  & -K_1^j\\
-K_2^j & 1+2K_2^j & -K_2^j &0& & 0\\
0&-K_3^j &1+2K_3^j&-K_3^j &&&\\
\vdots&&&\ddots&&\\
-K_N^j & \dots &&0&-K_N^j&1+2K_N^j
\end{bmatrix}
\begin{bmatrix}
\gamma_1^{j+1} \\ \\ \vdots \\ \\ \gamma_N^{j+1}
\end{bmatrix}
=
\begin{bmatrix}
\gamma_1^{j} \\ \\ \vdots \\ \\ \gamma_N^{j}
\end{bmatrix}.
\]

\subsection{The computation and results}\label{sec:computation-results}
To test the predictions of the $n$-loop conjecture, we approximated solutions to Curve Shortening Flow using the backward scheme \eqref{eq:linear-system}.  As for initial curve, we chose the curve parameterized by
\[
x=\sin \theta, \qquad y=\bigl(1-\lambda (2+\sin^2\theta)\bigr)\cos \theta
\]
varying the parameter $\lambda$ until we found a solution that comes close to shrinking to the origin while still remaining a $3$-loop.

We used Python, and, in particular, the NumPy \cite{harris2020array} and SciPy \cite{2020SciPy-NMeth} packages to implement the computations.

We discretized the curve to $N=5000$ points, and adjusted the time steps $\Delta t_j$ so that the constants $K_i^j$ never exceed $10^6$.  We stop the computation when the requirement that $\max K_i^j\leq 10^6$ forces us to choose $\Delta t_j< 10^{-9}$.

\begin{figure}[h]
\centering \includegraphics[width=\textwidth]{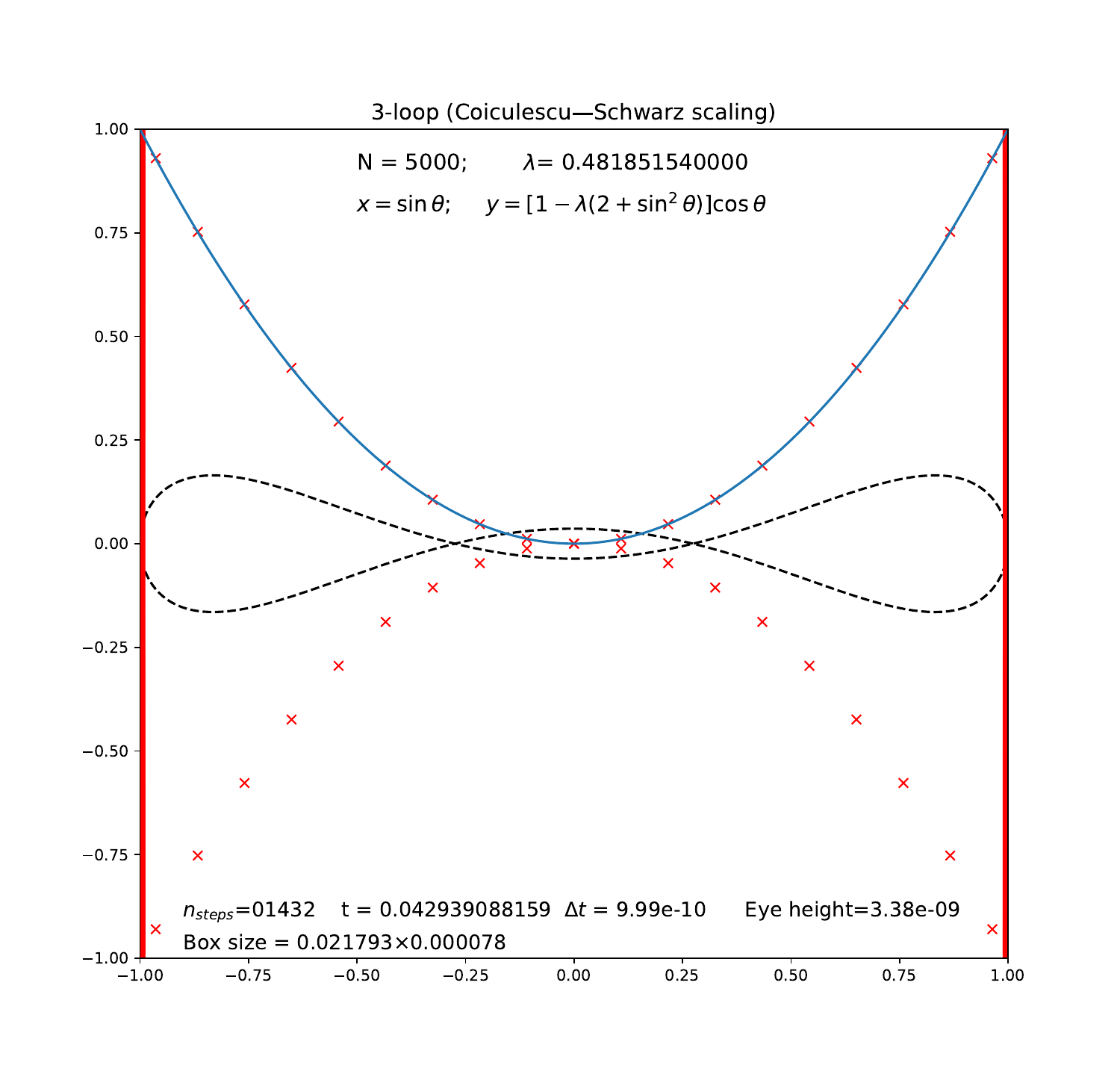} 
\caption{\textbf{Finding the shrinking 3-loop.}~~The solution right before the singular time $T_{\lambda_1}$ is represented by red crosses.  For comparison, the graph of $y=x^2$ is shown as a blue solid curve.  For these curves, the Coiculescu--Schwartz scaling is shown, i.e., ~they are scaled so as to fit the square.  In addition, the initial curves are drawn as black dashed curves in the standard Cartesian scaling. }
\label{fig:computed-3loop}
\end{figure}

\begin{figure}[h]
\includegraphics[width=\textwidth]{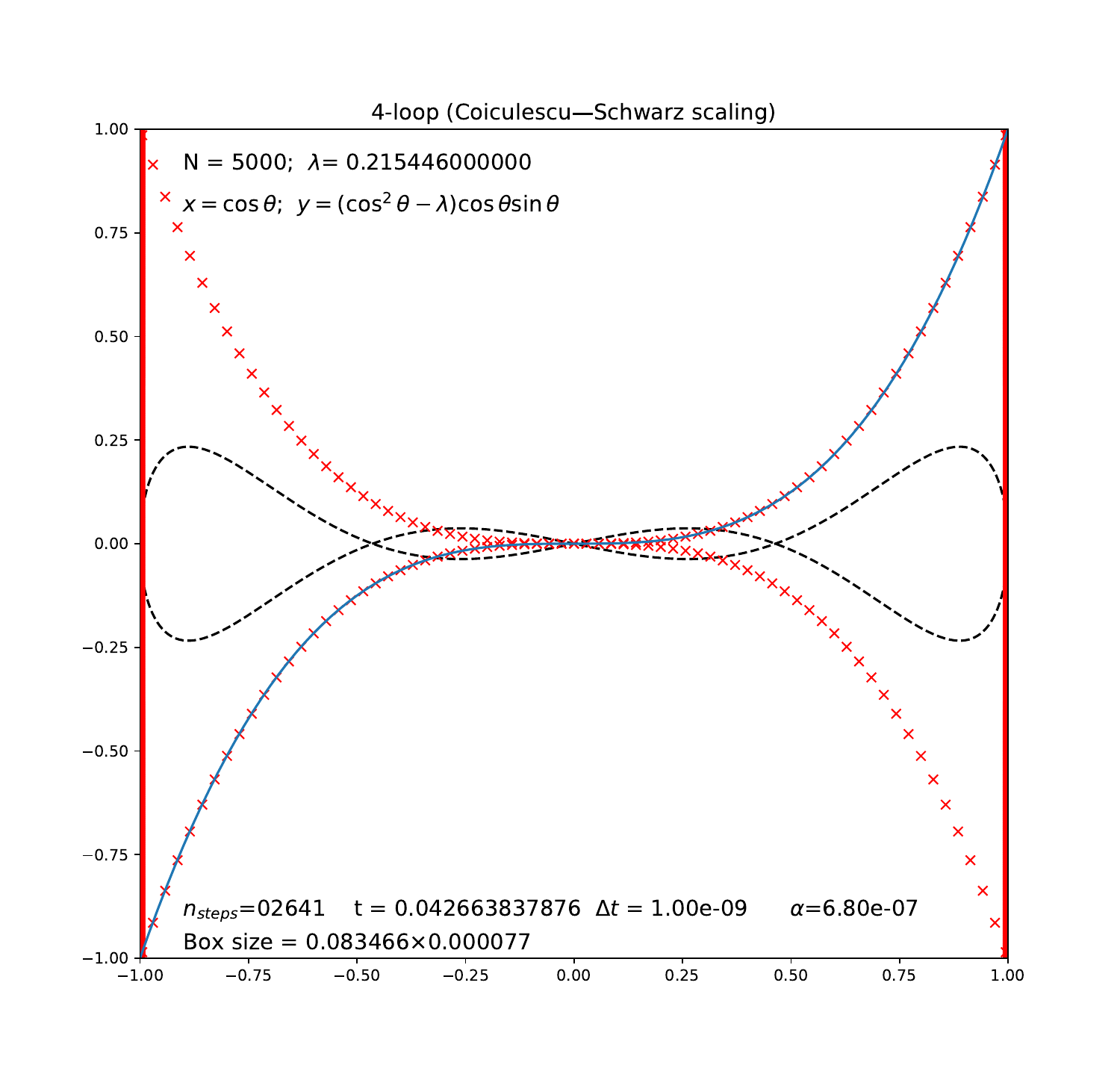}
\caption{\textbf{Finding the shrinking 4-loop.}~~As in Figure~\ref{fig:computed-3loop} we plot Coiculescu--Schwartz scalings of the solution right before it vanishes in red crosses, and the graph of $y=x^3$ as a solid blue curve.
 The initial curve is drawn as a black dashed curve in the standard Cartesian scaling.
}
\label{fig:computed-4loop}
\end{figure}

The choice $\lambda=0.48185154\dots$ led to a solution that shrinks to a bounding box of size $\approx 0.0218\times0.000078$. Then we plotted the solution in a square box (which provides the Coiculescu--Schwartz rescaling) and also drew the graph of $y=x^2$ in the same figure.  As one can see in Figures~\ref{fig:computed-3loop} and~\ref{fig:computed-4loop}, the computed solutions match the graphs of $y=x^2$ or $y=x^3$, respectively, and thus support the $n$-loop conjecture.

\bibliographystyle{plain}
\bibliography{refs}

\end{document}